\documentclass[12pt,letterpaper,twoside]{article}

\usepackage{amsmath,amssymb,amscd}
\usepackage{amsthm}
\usepackage{stmaryrd}
\usepackage{fancyhdr}
\usepackage{times}
\usepackage{mathrsfs} 
\usepackage{titlesec}
\usepackage[pdftex,dvips]{geometry}

\newtheoremstyle{plain}
     {\topsep}
     {\topsep}
     {\itshape}
     {}
     {\bfseries}
     {}
     { }
     {\thmnumber{#2.}\hspace{0.5ex}\thmname{#1.}\thmnote{ \rm (#3)}}

\newtheoremstyle{definition}
     {\topsep}
     {\topsep}
     {\normalfont}
     {}
     {\bfseries}
     {}
     { }
     {\thmnumber{#2.}\hspace{0.5ex}\thmname{#1.}\thmnote{ \rm (#3)}}

\newtheoremstyle{fact}
     {\topsep}
     {\topsep}
     {\slshape}
     {}
     {\bfseries}
     {}
     { }
     {\thmnumber{#2.}\hspace{0.5ex}\thmname{#1.}\thmnote{ \rm (#3)}}

\newtheorem{theorem}[subsection]{Theorem}
\newtheorem*{theorem*}{Theorem} 
\newtheorem{lemma}[subsection]{Lemma}
\newtheorem{proposition}[subsection]{Proposition}
\newtheorem{corollary}[subsection]{Corollary}

\theoremstyle{definition}

\newtheorem{remark}[subsection]{Remark}
\newtheorem*{remark*}{Remark}

\newtheorem*{question*}{Question}
\newtheorem*{examples*}{Examples}  

\newtheorem*{example*}{Example}

\newtheorem*{convention*}{Convention}

\theoremstyle{fact}

\def\proofont{\fontseries{bx}\fontshape{n}\selectfont}
\def\proofname{Proof.}

\newcommand{\pcite}[2]{{\cite[#1]{#2}}}

\newcommand{\rest}{\mbox{\parbox[t]{0.1cm}{$|$\\[-10pt] $|$}}}
\newcommand{\spec}{\operatorname{sp}}

\newcommand{\Note}[1]{}

\makeatletter
\renewenvironment{proof}[1][\proofname]{\par
  \normalfont
  \topsep6\p@\@plus6\p@ \trivlist
  \item[\hskip\labelsep\noindent\proofont #1]\ignorespaces
}{%
  \qed\endtrivlist
}
\makeatother

\titlelabel{\thetitle.\ }
\titleformat*{\section}{\normalsize\bfseries\centering}
\titleformat*{\subsection}{\normalsize\bfseries}
\titlespacing{\subsection}{0pt}{\topsep}{0.5ex}
\titleformat{\subsection}[runin]{\normalfont\bfseries}{%
\thesubsection.}{0.5ex}{}[.]

\author{Rachid El Harti and G\'abor Luk\'acs
\thanks{The second author gratefully acknowledges the generous financial support 
received from the Killam Trusts, and Dalhousie University that enabled him to 
do this research.}}
   
\title{Bounded and unitary elements in pro-$C^*$-algebras
}

\begin{document}

\makeatletter
\let\mytitle\@title
\chead{\small\itshape R. El Harti and G. Luk\'acs / \mytitle }
\fancyhead[RO,LE]{\small \thepage}
\makeatother

\maketitle

\def\thanks#1{} 

\thispagestyle{empty}

\nocite{Allan}
\nocite{Arveson}
\nocite{BhaKar}
\nocite{Board}
\nocite{Bourd2}
\nocite{Bourd1}
\nocite{EHR4}
\nocite{EHR1Chi}
\nocite{DubPor}
\nocite{Inoue}
\nocite{Koethe}
\nocite{Michael}
\nocite{Palmer2}
\nocite{Phil3}
\nocite{Phil2}
\nocite{Phil1}
\nocite{Schmu}
\nocite{Sebes}
\nocite{Zelaz}

\begin{abstract}
A {\em pro-$C^*$-algebra} is a (projective) limit of $C^*$-algebras 
in the category of topological $*$-algebras. From the perspective of 
non-commutative geometry, pro-$C^*$-algebras can be seen as non-commutative 
$k$-spaces.  An element of a pro-$C^*$-algebra is {\em bounded} if there
is a uniform bound for the norm of its images
under any continuous $*$-homomorphism into a $C^*$-algebra. 
The $*$-subalgebra consisting of the bounded elements turns out
to be a $C^*$-algebra.  In this paper, we investigate pro-$C^*$-algebras from 
a categorical point of view. We study the functor $(-)_b$ that assigns to 
a pro-$C^*$-algebra the $C^*$-algebra of its bounded elements, which is
the dual of the Stone-\v Cech-compactification.
We show that $(-)_b$ is a coreflector, and it preserves exact sequences. 
A generalization of the Gelfand duality for
commutative unital pro-$C^*$-algebras is also presented.
\end{abstract}

\section{Introduction}

\label{sect:intro}

In a $C^*$-algebra, the norm is determined by the algebraic structure, and 
thus being a $C^*$-algebra is an entirely algebraic property. This 
intuitive conclusion was confirmed by van Osdol, who proved that the unit 
ball functor, sending a $C^*$-algebra $A$ to the {\em set}
$\{a \in A \mid \|a| \leq 1\}$, is monadic (cf. \cite{Osdol}). 
Later on, Pelletier and  Rosick{\'y} investigated the
system of  operations and equations representing the unit balls of
$C^*$-algebras as universal algebras (cf. \cite{PelRos} and \cite{PelRos2}).

A Hausdorff {\em $k$-space} is a colimit of compact Hausdorff spaces in 
the category of Hausdorff spaces (and their continuous maps).  
Such spaces were thoroughly investigated by Brown, Steenrod, and Dubuc and 
Porta (cf. \cite[3.3]{Brown}, \cite{Steenrod}, and \cite{DubPor}).
A {\em pro-$C^*$-algebra} is a (projective) limit of 
$C^*$-algebras in the category of topological $*$-algebras.
Such algebras were studied under various names ($LMC^*$-algebras, locally 
$C^*$-algebras, and $\sigma$-$C^*$-algebra in 
the metrizable case) by Schm\"udgen \cite{Schmu}, Inoue \cite{Inoue},
Arveson \cite{Arveson}, Phillips \cite{Phil2}, and
El Harti \cite{EHR4}. 

By the Gelfand duality, the category of commutative unital $C^*$-algebras 
is equivalent to the opposite of the category of compact Hausdorff spaces. 
It turns out that using a close relative of $k$-spaces, the Gelfand 
duality can be extended to commutative unital pro-$C^*$-algebras 
(Theorem~\ref{thm:genGel}). Thus, from the perspective of non-commutative 
geometry, pro-$C^*$-algebras can be seen as non-commutative $k$-spaces.

\bigskip

An element $a$ of a pro-$C^*$-algebra $A$ is {\em bounded} if
there is a constant $M$ such that $\|\varphi(a)\|\leq M$ for 
every continuous $*$-homomorphism $\varphi\colon A \rightarrow B$ 
into a $C^*$-algebra $B$. The smallest constant $M$ with 
this property is denoted by $\|a\|_\infty$, and the 
$*$-subalgebra $A_b$ of bounded elements in $A$ is a
$C^*$-algebra with respect to $\|\cdot\|_\infty$ (see
section~\ref{sect:b}). Since the Stone-\v Cech-compactification
of a Tychonoff space $X$ can be obtained as the space
of multiplicative functionals of the $C^*$-algebra
$C_b(X)$ with respect to the norm of uniform convergence, 
it is natural to think of the functor $A \mapsto A_b$ 
as the non-commutative dual of the Stone-\v Cech-compactification
functor $X \mapsto \beta X$.

\bigskip

In this paper, we investigate pro-$C^*$-algebras from a categorical
point of view. After a brief presentation of terminology and basic
results on pro-$C^*$-algebras in section~\ref{sect:preli}, we proceed
by presenting a generalization of the Gelfand duality for
commutative unital pro-$C^*$-algebras in section~\ref{sect:ab}.
We also compare this generalized Gelfand duality with results of 
Dubuc and Porta \cite{DubPor}.
In section~\ref{sect:b}, we study the functor $(-)_b$ that assigns
to a pro-$C^*$-algebra $A$ the $C^*$-algebra $A_b$ of its bounded elements. 
We show that this functor is a coreflector (Theorem~\ref{thm:b:coref}), 
and it preserves exactness of sequences (Theorem~\ref{thm:Ab:ex}). 
Finally, in section~\ref{sect:u}, we make a few observations
concerning the connected component of the identity of the group
of unitary elements in a pro-$C^*$-algebra.

\section{Preliminaries}

\label{sect:preli}

In this section, we introduce some terminology, and present well-known 
elementary results on pro-$C^*$-algebras.

\subsection{Topological $\boldsymbol *$-algebras} 
A {\em $*$-algebra} (or {\em involutive algebra}) 
is an algebra $A$ over $\mathbb{C}$ with an
involution $\empty^*\colon A \rightarrow A$ such that
$(a+\lambda b)^* = a^*+\bar \lambda b^*$ and $(ab)^*=b^*a^*$
for every $a,b\in A$ and $\lambda \in \mathbb{C}$.
A linear seminorm 
$p$ on a $*$-algebra $A$ is a {\em $C^*$-seminorm} if
$p(ab) \leq p(a) p(b)$ ({\em submultiplicative}) and 
$p(a^* a)=p(a)^2$ ({\em $C^*$-condition}) for every $a,b\in A$. 
Note that the $C^*$-condition alone implies that $p$ is 
submultiplicative, and in particular 
$p(a^*)=p(a)$ for every $a \in A$ (cf. \cite{Sebes}, \cite[9.5.14]{Palmer2}).
A {\em topological $*$-algebra} is a $*$-algebra $A$ equipped with
a topology making the operations (addition, multiplication,
additive inverse, involution) jointly continuous.
The category of topological $*$-algebras
and their continuous $*$-homomorphisms is denoted by $\mathsf{T^*A}$.
For $A \in \mathsf{T^*A}$, one puts $\mathcal{N}(A)$ for the set of 
continuous $C^*$-seminorms on $A$;  $\mathcal{N}(A)$ is a directed set
with respect to pointwise ordering, because 
$\max\{p,q\}\in\mathcal{N}(A)$ for every $p,q\in\mathcal{N}(A)$.

\subsection{$\boldsymbol {C^*}$-algebras}
A {\em $C^*$-algebra} is a complete Hausdorff topological algebra
whose topology is given by a single $C^*$-norm. The full subcategory of
$\mathsf{T^*A}$ formed by the $C^*$-algebras is denoted by $\mathsf{C^*A}$.
For $A\in \mathsf{T^*A}$ and $p\in \mathcal{N}(A)$,
$\ker p = \{ a\in A \mid p(a)=0\}$ is a $*$-ideal in $A$, and
$p$ induces a $C^*$-norm on the quotient $A/\ker p$, so
the completion $A_p$ of this quotient with respect to
$p$ is a $C^*$-algebra. Each pair $p,q \in \mathcal{N}$ such that
$q \geq p$ induces a natural (continuous) surjective $*$-homomorphism
$\pi_{pq}\colon A_q \rightarrow A_p$, which turns 
$A_{(-)} \colon \mathcal{N}(A) \longrightarrow \mathsf{C^*A}$ into a functor.

\subsection{The $\boldsymbol *$-representation topology}
For $A \in \nolinebreak\mathsf{T^*A}$, a {\em $*$-representation} of $A$
is a continuous $*$-homomorphism $\pi\colon A\rightarrow B(H)$
of $A$ into the $C^*$-algebra of bounded operators on some Hilbert space $H$
(i.e., a morphism in $\mathsf{T^*A}$). The class
of $*$-representations of $A$ is denoted by $\mathcal{R}(A)$.
By the Gelfand-Na\u{\i}mark-Segal theorem,
every $C^*$-algebra is $*$-isomorphic (and thus isometric) to a closed
subalgebra of $B(H)$ for a large enough Hilbert space $H$
(cf. \cite[2.6.1]{DixmierC}), so each $A_p$ embeds (isometrically) into $B(H_p)$
for some Hilbert space $H_p$. Thus, for each $p\in \mathcal{N}(A)$ one obtains 
a $*$-representation
$\bar\pi_p\colon A \rightarrow A/\ker p \rightarrow A_p \rightarrow B(H_p)$
such that $p(x)=\|\bar\pi_p(x)\|$. Conversely, each $\pi \in \mathcal{R}(A)$
gives rise to a $C^*$-seminorm $p_\pi(x)=\|\pi(x)\|$. Therefore, the initial
topology $\mathcal{T}_A$ induced by the class $\mathcal{R}(A)$ coincides
with the one induces by the family of $C^*$-seminorms $\mathcal{N}(A)$.
The topology $\mathcal{T}_A$ is called the {\em $*$-representation topology}.

\subsection{Pro-$\boldsymbol {C^*}$-algebras}
A {\em pro-$C^*$-algebra} is a complete Hausdorff topological $*$-algebra 
$A$ whose topology is defined by a family of $C^*$-seminorms, or 
equivalently, the topology of $A$ coincides with $\mathcal{T}_A$, and 
$\mathcal{T}_A$ is Hausdorff and complete (i.e., every Cauchy filter or 
net converges). For $A\in\mathsf{T^*A}$, the following statements are 
equivalent (cf. \cite[1.1.1]{Phil1}):

\begin{list}{{\rm (\roman{enumi})}}
{\usecounter{enumi}\setlength{\labelwidth}{25pt}\setlength{\topsep}{-8pt}
\setlength{\itemsep}{-5pt} \setlength{\leftmargin}{25pt}}

\item
$A \cong \lim\limits_{\stackrel {\longleftarrow\!\!-\!\!-} 
{\empty_{p\in \mathcal{N}(A)}}}  A_p$ in $\mathsf{T^*A}$;

\item
$A$ is a (projective) limit of $C^*$-algebras in $\mathsf{T^*A}$;

\item
$A$ is a closed $*$-subalgebra of a product of $C^*$-algebras 
in $\mathsf{T^*A}$;

\item
$A$ is a pro-$C^*$-algebra.
\end{list}
\vspace{5.5pt}
It turns out that for a pro-$C^*$-algebra $A$, $(A/\ker p, p)$ is already
complete, so $A_p$ is a quotient of $A$ (cf. \cite[Thm.~2.4]{Apost},
\cite[Folg.~5.4]{Schmu}, \cite[1.2.8]{Phil1}, \cite[1.12]{Phil2}). 
Pro-$C^*$-algebras form a full reflective subcategory of $\mathsf{T^*A}$ 
that we denote by $\overline{\mathsf{P^*A}}$ (cf. \cite[5.1]{GL7}).
In light of (i), we often refer to an element $a$ of a pro-$C^*$-algebra 
$A$ as $(a_p)_{p \in \mathcal{N}(A)}$, where $a_p \in A_p$ and 
$\pi_{pq}(a_q)=a_p$ for every $p,q \in \mathcal{N}(A)$.

\subsection{Examples} \label{exm:basics}

\begin{list}{{\rm (\Alph{enumi})}}
{\usecounter{enumi}\setlength{\labelwidth}{25pt}\setlength{\topsep}{-8pt}
\setlength{\itemsep}{-5pt} \setlength{\leftmargin}{25pt}}

\item
$\prod\limits_{n=1}^\infty \mathbb{M}_n(\mathbb{C})$, the
product of the (full) matrix algebras  $\mathbb{M}_n(\mathbb{C})$ in 
$\mathsf{T^*A}$, is a metrizable pro-$C^*$-algebra.

\item
The algebra $C_\omega([0,1])$ of continuous complex-valued maps on $[0,1]$
with the topology of uniform convergence on compact countable subsets
is a pro-$C^*$-algebra (for details, see section~\ref{sect:ab}).

\item
The tangent algebra of a $C^*$ introduced by Arveson
(defined as the universal factorizer of derivations into $C^*$-algebras)
fails to be $C^*$-algebra, but it is a metrizable pro-$C^*$-algebra
(for details, see \cite[5.9]{Arveson} and \cite[2.1]{Phil1}).

\end{list}
\vspace{5.5pt}

\subsection{Spectrum and unitization}
We denote by $G_A$ the set of invertible elements in an algebra $A$.
The {\em spectrum} of an element $a$ in a unital algebra $A$ over 
$\mathbb{C}$ is defined as 
\begin{align} 
\spec_A(a) &=\{ \lambda \in \mathbb{C} \mid 
a - \lambda 1  \not \in G_A \}.
\intertext{If $A \in \overline{\mathsf{P^*A}}$, then $a=(a_p) \in A$ is 
invertible if and only if each $a_p$ is invertible in $A_p$. Therefore,}
\spec_A(a) & = \hspace{-8pt}\bigcup\limits_{p\in\mathcal{N}(A)} 
\hspace{-8pt} \spec_{A_p}(a_p). \label{eq:specPA}
\end{align}
If $B \subseteq A$ is a closed $*$-subalgebra that contains the unit of $A$ and 
$b \in B$, then $\spec_B(b)=\spec_A(b)$ (cf. \cite[1.2]{Schmu}).
If $A$ is not unital (which is the case for ideals in a pro-$C^*$-algebra), 
one can adjoin a unit by putting
$A^+=A\oplus \mathbb{C}e$, extending multiplication in the obvious way,
and setting $(a,\lambda)^*=(a^*,\bar \lambda)$. (Certainly, 
$A^+ = \lim\limits_{\longleftarrow} A_p^+$ is a pro-$C^*$-algebra.)
In this case, one defines $\spec_{A^+}(a)$ as the spectrum of $a$. 
This definition is consistent (up to $\{0\}$), because if $A$
already had a unit, then $\spec_{A^+} (a) = \spec_A(a)\cup \{0\}$.
Thus, if $I$ is a closed (two-sided) $*$-ideal of a unital 
pro-$C^*$-algebra $A$ and $c \in I$,
then $\spec_{I^+}(c)=\spec_{A^+}(c)=\spec_A (c)$, because $c$
cannot be invertible in $A$.

\subsection{Functional calculus}
Let $A$ be a unital pro-$C^*$-algebra and let $a =(a_p)\in A$.
Suppose that $f$ is a continuous complex-valued map
on $\spec_A (a)$ and either $a$ is {\em normal} (i.e., 
$a^*a =a a^*$) or $f$ is analytic on a neighborhood
of $\spec_A (a)$. It follows from (\ref{eq:specPA}) that
$f$ is continuous on each $\spec_{A_p}(a_p)$, and either
each $a_p$ is normal or $f$ is analytic on a neighborhood of each 
$\spec_{A_p}(a_p)$. In both cases, $f(a_p)$ is defined,
$\spec_{A_p}(f(a_p))=f(\spec_{A_p}(a_p))$, and 
$\pi_{pq}(f(a_q))=f(a_p)$ for each pair $q \geq p$ in 
$\mathcal{N}(A)$. Therefore, $(f(a_p))\in A$ and
$\spec_A(f(a))=f(\spec_A(a))$. In case $A$ is not unital, 
$f(0)=0$ is the only additional assumption needed in 
order to ensure that $f(a)$ (calculated in $A^+$) belongs 
to~$A$.

\subsection{Automatic continuity fails in $\boldsymbol{\overline{\mathsf{P^*A}}}$}
It is a well-known fact that every $*$-ho\-mo\-mor\-phism of an involutive 
Banach algebra into a $C^*$-algebra is continuous (furthermore, its norm is 
$\leq 1$; cf. \cite[1.3.7]{DixmierC}). In particular, every
$*$-homomorphism between $C^*$-algebras is continuous. Results of this nature 
are referred to as ``automatic continuity" of homomorphisms, because the 
continuity is deduced from a purely algebraic assumption 
(cf. \cite[Chapter 6]{Palmer1}). Unfortunately, in general, this result
fails for pro-$C^*$-algebras: Both $C(\omega_1)$ and
$C(\omega_1 +1)$ are pro-$C^*$-algebras in the topology of uniform 
convergence on compacta, but the  $*$-homomorphism
$C(\omega_1)\rightarrow C(\omega_1 +1)$ given by setting
$f(\omega_1)=\lim\limits_{x\rightarrow\omega_1} f(x)$ fails to be
continuous (cf. \cite[12.2]{Michael}, \cite[1.4.9]{Phil1}).

\subsection{Automatic continuity for 
$\boldsymbol \sigma$-$\boldsymbol C^*$-algebras} \label{sub:auto:sigma}
A {\em $\sigma$-$C^*$-algebra} is a pro-$C^*$-algebra whose topology
is metrizable, or equivalently, it is a (projective) limit of a sequence of 
$C^*$-algebras in $\mathsf{T^*A}$. Every $*$-homomorphism of 
a $\sigma$-$C^*$-algebra into a pro-$C^*$-algebra is continuous
(cf. \cite[5.2]{Phil2}, \cite[1.1.6]{Phil1}).

\section{Commutative pro-$\boldsymbol{C^*}$-algebras}

\label{sect:ab}

By the Gelfand duality, every commutative unital $C^*$-algebra
is $*$-isomorphic (and thus isometrically isomorphic) to the algebra
of continuous maps on the set of its multiplicative functionals (i.e.,
characters). In this section, we present some generalizations of this result
to pro-$C^*$-algebras. Our aim is to adhere to a categorical point of view,
and thus it slightly differs from the approach of Inoue \cite{Inoue} and
Phillips \cite[1.4]{Phil1} \& \cite[2]{Phil2}, who obtained essentially the 
same results.

\subsection{Leading example}
For a Hausdorff space $X$, one puts $\mathcal{K}(X)$ for the collection
of compact subsets of $X$. For $\mathcal F \subseteq \mathcal K (X)$ such
that $\bigcup \mathcal{F} = X$, we say that a map $f\colon X \rightarrow Y$
is $\mathcal{F}$-continuous if $f\rest_F$ is continuous for every 
$F\in \mathcal{F}$. We put $C_\mathcal{F}(X)$ for the $*$-algebra of complex-valued 
$\mathcal{F}$-continuous maps $f$ on $X$, and provide it with the
topology of {\em uniform convergence on $\mathcal{F}$} given
by the family of $C^*$-seminorms $\{p_F\}_{F\in \mathcal F}$, where
$p_F(f) = \sup\limits_{x\in F} |f(x)|$. It is a well-known result on
uniform spaces that $C_\mathcal{F}(X)$ is complete (cf. \cite[7.10(d)]{Kelley}),
and thus $C_\mathcal{F}(X)$ is a pro-$C^*$-algebra. We say that 
$\mathcal{F}$ is a {\em strongly functionally generating} family for $X$ if
every map in $C_\mathcal{F}(X)$ is continuous (cf. \cite[p.~111]{ArhBook}).
The Gelfand duality states that every commutative unital $C^*$-algebra
is isomorphic in $\mathsf{C^*A}$ (and thus in $\mathsf{T^*A}$) to the 
algebra $C(K)$ of continuous maps with the topology of uniform 
convergence on some compact Hausdorff space 
$K$. We show that every commutative unital pro-$C^*$-algebra is isomorphic 
to $C_\mathcal{F}(X)$ in $\mathsf{T^*A}$ for some space $X$, 
where $X$ is strongly functionally generated by $\mathcal{F}$.

\subsection{The categories $\boldsymbol{\mathsf{kHaus}}$ and 
$\boldsymbol{\mathsf{k_RTych}}$} \label{sub:k}
A Hausdorff space $X$ is a {\em $k$-space} if every $\mathcal{K}(X)$-continuous 
(i.e., {\em $k$-continuous}) map 
$f\colon X \rightarrow Y$ into a Hausdorff space $Y$ is continuous, or 
equivalently, if $U \subseteq X$ is open provided that $U\cap K$ is open 
in $K$ for every compact subset $K$ of $X$ (cf. \pcite{3.3.21}{Engel6}). 
The category $\mathsf{kHaus}$ of Hausdorff $k$-spaces and continuous maps
is cartesian closed, and it is a coreflective subcategory of the category
$\mathsf{Haus}$ of Hausdorff spaces (cf. \cite[3.3]{Brown}, 
\cite[VII.8]{MacLane}). (The coreflector $\mathsf{k}\colon 
\mathsf{Haus}\longrightarrow \mathsf{kHaus}$ is often referred to as 
{\em $k$-ification}. Some authors use the term ``compactly generated 
spaces" for $k$-spaces.)
A Hausdorff space $X$ is a {\em $k_R$-space} if it is strongly
functionally generated by $\mathcal{K}(X)$ (i.e., every
real-valued $k$-continuous map $f\colon X \rightarrow \mathbb{R}$ on $X$ is 
continuous). Luk\'acs showed that the category $\mathsf{k_RTych}$ of Tychonoff
$k_R$-spaces and continuous maps is cartesian closed, $\mathsf{k_RTych}$ 
is a coreflective subcategory of $\mathsf{Tych}$ (Tychonoff spaces), and 
$\mathsf{k_RTych}$ is equivalent to a full epireflective subcategory of 
$\mathsf{kHaus}$  (cf. \cite[Theorems 3.2, 2.1, 3.1]{GL5}). The motivation 
for introducing  this category is the incompatibility of the $k$-space 
property  and the Tychonoff property (cf. \cite[Example]{GL5}).  Their 
necessity for pro-$C^*$-algebras roots in the observation that if $X$ 
admits a strongly functionally generating family of compact subsets, then 
it must be a $k_R$-space.

\subsection{The character space $\boldsymbol{\Delta(A)}$} \label{sub:Delta}
For a topological vector space $X$, one puts $X^\prime$ for the {\em dual space}
of $X$ consisting of the continuous linear functionals on $X$. The
{\em $w^*$-topology} on $X^\prime$ is the topology of pointwise convergence
on the points of $X$---in other words, $X^\prime$ is provided with the initial
topology with respect to all evaluations maps 
$\{\operatorname{ev}_x\colon X^\prime \rightarrow \mathbb{C}\}_{x\in X}$.
Since $(X^\prime,w^*)$ is a Hausdorff topological vector space, 
it is Tychonoff. 
For $A \in \mathsf{T^*A}$, one puts $\Delta(A)$ for the subspace of 
$(A^\prime,w^*)$ consisting of the non-zero $*$-homomorphisms. 
Using this notation, Gelfand duality states that the map 
$\operatorname{ev}_A \colon A \rightarrow C(\Delta(A))$, defined by $a 
\mapsto \operatorname{ev}_a$, is an isomorphism of $C^*$-algebras for 
every  $C^*$-algebra $A$.

\begin{theorem} \label{thm:genGelf:1half}
Let $A$ be a commutative unital pro-$C^*$-algebra. Then:

\begin{list}{{\rm (\alph{enumi})}}
{\usecounter{enumi}\setlength{\labelwidth}{25pt}\setlength{\topsep}{-8pt}
\setlength{\itemsep}{-5pt} \setlength{\leftmargin}{25pt}}

\item
$\Delta(A)=\hspace{-10pt}\bigcup\limits_{p \in \mathcal{N}(A)} \hspace{-10pt}
\Delta (A_p)$;

\item
$\Delta(A)$ is a Tychonoff $k_R$-space, and 
$\Phi(A)=\{\Delta(A_p) \mid p \in \mathcal{N}(A)\}$ is a strongly 
functionally generating family;

\item
the natural map $\operatorname{ev}_A \colon A \rightarrow
C_{\Phi(A)}(\Delta(A))$ defined by $a \mapsto \operatorname{ev}_a$ is an 
isomorphism of pro-$*$-algebras.

\end{list}

\end{theorem}

\begin{proof} (a) By \ref{sub:Delta}, $\Delta(A)$ is a subspace of the 
Tychonoff space $A^\prime$, and so it is Tychonoff. 
Since $A$ is commutative and unital,
so are its quotients $A_p$. Thus, $A_p\cong C(\Delta(A_p))$ by Gelfand duality, 
and the $\Delta(A_p)$ are compact.  Each project $A \rightarrow A_p$
gives rise to a continuous injective map $\Delta(A_p) \rightarrow \Delta(A)$,
which is an embedding, because $\Delta(A_p)$ is compact. On the other hand,
if $\chi \in \Delta(A)$, then $p_{\chi}(a)=|\chi(a)|$ is a continuous 
$C^*$-seminorm on $A$, and therefore $\chi \in \Delta(A_{p_\chi})$.

(b) Suppose that $f\colon \Delta(A) \rightarrow \mathbb{R}$ is 
$\Phi(A)$-continuous, in other words, $f\rest_{\Delta(A_p)}$ is continuous
for every $p \in \mathcal{N}(A)$. Then, by Gelfand duality, for each
$p \in \mathcal{N}(A)$ there is $a_p \in A_p$ such that
$f\rest_{\Delta(A_p)}\!\! = \operatorname{ev}_{a_p}$. Furthermore, 
since $f\rest_{\Delta(A_q)} \rest_{\Delta(A_p)} \!\! = 
f\rest_{\Delta(A_p)}\!\!$, one has $\pi_{pq}(a_q)=a_p$
for every $q \geq p$. Thus, $a=(a_p) \in 
A$ and $f=\operatorname{ev}_a$.
In particular, $f$ is continuous, and therefore $\Phi(A)$ is a strongly 
functionally generating family. Hence, $\Delta(A)$ is a $k_R$-space.

(c) It is clear that $\operatorname{ev}_A$ is a $*$-homomorphism. Since
\begin{equation}
p(a)=\sup \{ |\rho(a)| \mid \rho \in \Delta(A_p)\} = 
p_{\Delta(A_p)}(\operatorname{ev}_a)
\end{equation}
for every $p \in \mathcal{N}(A)$, it follows that $\operatorname{ev}_A$ 
is continuous and open onto its image.
If $\operatorname{ev}_A(a)=0$ for $a \in A$, then $\rho(a)=0$ for every 
$\rho \in \Delta(A)$, and so $p(a)=0$ for every $p \in \mathcal{N}(A)$.
Thus, $a=0$ (because $A$ is Hausdorff). If $f \in C_{\Phi(A)}(\Delta(A))$, 
then by the argument developed in (b), $f=\operatorname{ev}_a$ for some 
$a \in A$. Therefore, $\operatorname{ev}_A$ is surjective. Hence, 
$\operatorname{ev}_A$ is an isomorphism of pro-$C^*$-algebras.
\end{proof}

\begin{lemma} \label{lemma:genGel:2half}
Let $X$ be a Tychonoff $k_R$-space, and let $\mathcal{F}$ be a strongly 
functionally generating family of compact subsets of $X$. Then the natural 
map $\operatorname{ev}_X\colon X \rightarrow \Delta(C_\mathcal{F} (X))$ 
defined by $x \mapsto \operatorname{ev}_x$ is a homeomorphism.
\end{lemma}

\begin{proof}
Since $X$ is a Tychonoff space $k_R$-space, it embeds into 
$C_\mathcal{F}(X)^\prime$ by evaluation (cf. \cite[Lemma~2.1]{GL5}), and 
its image is contained in
$\Delta(C_\mathcal{F} (X))$, so $\operatorname{ev}_X$ is an embedding.
In order to show that $\operatorname{ev}_X$ is surjective, let 
$\rho \in \Delta(C_\mathcal{F}(X))$. Since $\rho$ is a continuous 
functional with respect to the topology determined by the seminorms
$\{p_F\}_{F\in \mathcal{F}}$,  there are
$F_1,\ldots, F_l \in \mathcal{F}$ and a constant $r$ such that 
$|\rho(f)| \leq r \max\{p_{F_1}(f),\ldots,p_{F_l}(f)\}$ for every 
$f \in C_\mathcal{F}(X)$. Thus, for $F_0=F_1 \cup \ldots \cup F_l$, 
$\rho$ factors through a multiplicative functional $\bar 
\rho$ of the $C^*$-algebra $C(F_0)$, that is, 
$\rho(f)=\bar \rho (f \rest_{F_0})$. Therefore, by 
Gelfand duality, there is $x_0 \in F_0$ such that 
$\bar\rho(g)=g(x_0)$ for 
every $g \in C(F_0)$, and hence $\rho=\operatorname{ev}_X(x_0)$, as 
desired.
\end{proof}

\subsection{Notations} The category $\mathsf{SFG}$ is defined as follows:
Objects are pairs $(X,\mathcal{F})$, where $X$ is a Tychonoff $k_R$-space 
and $\mathcal{F}$ is a strongly functionally generating family of compact 
subsets of $X$. A morphism $g\colon (X,\mathcal{F})\rightarrow 
(Y,\mathcal{G})$ in $\mathsf{SFG}$ is a continuous map with the property 
that for every $F \in \mathcal{F}$, there is a finite subfamily	
$\{G_1,\ldots,G_l\}\subseteq \mathcal{G}$ such that 
$g(F)\subseteq G_1 \cup \ldots G_l$. We denote by 
$\mathsf{C_1 \overline{P^*A}}$ the full subcategory of 
$\mathsf{\overline{P^*A}}$ formed by the commutative unital 
pro-$C^*$-algebras.

\begin{theorem}  \label{thm:genGel}
The pair of functors
\begin{align}
\mathsf{C_1 \overline{P^*A}}^{\mathrm{op}} & 
\longrightarrow \mathsf{SFG} \\
A & \longmapsto (\Delta(A),\Phi(A))\\
C_\mathcal{F}(X) & \longmapsfrom (X,\mathcal{F})
\end{align}
form an equivalence of categories.
\end{theorem}

\begin{proof}
By Theorem~\ref{thm:genGelf:1half}(c), $A \cong C_{\Phi(A)}(\Delta(A))$ 
for every $A \in \mathsf{C_1 \overline{P^*A}}$. On the other hand, by 
Lemma~\ref{lemma:genGel:2half}, $X \cong \Delta(C_\mathcal{F}(X))$ 
for every $(X,\mathcal{F}) \in \mathsf{SFG}$. Furthermore, 
$F \in \Phi(C_\mathcal{F}(X))$ for every $F \in \mathcal{F}$,
because $p_F \in \mathcal{N}(C_\mathcal{F}(X))$. In order to complete the 
proof, let $G \in \Phi(C_\mathcal{F}(X))$ (and we show that it is 
contained in a finite union of members of $\mathcal{F}$). Then 
$p_G \in \mathcal{N}(C_\mathcal{F}(X))$, and thus there are 
$F_1,\ldots, F_l \in \mathcal{F}$ and a constant $r$
such that $p_G \leq r \max\{p_{F_1},\ldots,p_{F_l}\}$, 
because $\{p_F\}_{F \in \mathcal{F}}$ generates 
the topology of $C_\mathcal{F}(X)$. Since $X$ is Tychonoff,
if $x_0 \not\in F_1 \cup \ldots \cup F_l$, then there exists a continuous 
map $h$ on $X$ such that $h(x_0)=1$ and $h\rest_{F_i}= 0$ for each
$i$. So, $p_{F_i}(h)=0$, and therefore $p_G(h)=0$. In particular, 
$x_0 \not\in G$, and hence $G \subseteq F_1 \cup\ldots\cup F_l$, as 
desired.
\end{proof}

\subsection{$\mathsf{k_R}$-ification using function spaces}
Recall that for a Tychonoff space $X$, the Stone-\v Cech-compactification 
$\beta X$ (i.e., the reflection of $X$ into the category of compact 
Hausdorff spaces) can be computed using Gelfand duality:  
$\beta X \cong \Delta(C_b(X))$, where 
$C_b(X)$ is the $C^*$-algebra of all continuous bounded maps on $X$ with 
the norm of uniform convergence. Theorem~\ref{thm:genGel} provides 
a similar method for finding the $k_R$-ification
$\mathsf{k_R}X$ of a Tychonoff space $X$ (that is, the coreflection of $X$
into $\mathsf{k_RTych}$): By \cite[Lemma~2.1(b)]{GL5}, 
$\mathsf{k_R}X\cong \Delta(C_{\mathcal{K}(X)}(X))$. Naturally, 
when $X$ is not a $k_R$-space, $C_{\mathcal{K}(X)}(X)$ contains 
some non-continuous functions that are only $k$-continuous on $X$.

\subsection{$\boldsymbol K$-algebras} Dubuc and Porta 
investigated the category of $*$-algebra objects in $\mathsf{kHaus}$, 
called $K$-algebras (cf. \cite{DubPor}). 
These are $*$-algebras on a Hausdorff $k$-space $A$, with continuous 
unitary operations, and with binary operations being continuous on $A 
\times_\mathsf{k} A =\mathsf{k} (A \times A)$ (which is the cartesian 
product in $\mathsf{kHaus}$). So, addition and multiplication are expected 
to be only $k$-continuous on $A \times A$. A typical example of a 
$K$-algebra is $\mathsf{k}A$, the $k$-ification of a  pro-$C^*$-algebra $A$.
Dubuc and Porta established a 
dual adjunction between commutative unital $K$-algebras  and {\em functionally 
Hausdorff} $k$-spaces (i.e., continuous real-valued maps separate points) 
using the cotensor of the $\mathsf{kHaus}$-enriched category of 
$K$-algebras, and showed that the  Gelfand duality is a restriction of 
this dual adjunction. They also proved that the $k$-ification of a 
commutative pro-$C^*$-algebra admits a canonical monad-algebra structure 
with respect to the resulting monad (cf. \cite[3.8]{DubPor}). An essential 
difference between this setting and the one present in this paper is that 
in the context of $K$-algebras, the natural topology on $\Delta(A)$ is the 
$k$-ified compact-open topology, and not the $w^*$-topology.

Our approach is akin to the methods employed by Dubuc and Porta in a later 
paper, and by Phillips, using {\em quasitopological spaces} (cf. 
\cite{DubPor2} and \cite{Phil2}). For instance,  Phillips showed that 
$\mathsf{C_1 \overline{P^*A}}^{\mathrm{op}}$ is equivalent to the category 
$\mathsf{QHaus_f}$ of functionally Hausdorff quasitopological spaces 
(cf.~\cite[2.7]{Phil2}). Thus, by Theorem~\ref{thm:genGel}, we conclude:

\begin{corollary}
The categories $\mathsf{SFG}$ and $\mathsf{QHaus_f}$ are equivalent.
\qed
\end{corollary}

\section{The algebra of bounded elements}

\label{sect:b}

For a Tychonoff $k_R$-space $X$, the $C^*$-algebra $C_b(X)$ of continuous 
complex-valued bounded maps on $X$ is precisely the subalgebra of 
$C(X)$ of the maps $f\colon X \rightarrow \mathbb{C}$ such that
$\hspace{-6pt}\sup\limits_{K \in \mathcal{K}(X)}\hspace{-6pt} p_K(f) < \infty$. 
In this section, a natural generalization  of this notion 
to the non-commutative case is presented.

\subsection{Bounded elements} \label{sub:Ab}
For a pro-$C^*$-algebra $A$, 
$\|a\|_\infty\hspace{-2pt}=\hspace{-6pt}\sup\limits_{p \in \mathcal{N}(A)} 
\hspace{-6pt} p(a)$ is called the {\em uniform norm}, and $a \in A$ is
{\em bounded} if $\|a\|_\infty < \infty$. The $*$-subalgebra $A_b$ of 
the bounded elements in $A$ is dense in $A$, and $A_b$ is is a 
$C^*$-algebra  with the norm  $\|\cdot\|_\infty$ 
(cf. \cite[Satz~3.1]{Schmu}, \cite[1.2.4, 1.2.7]{Phil1},
\cite[1.11]{Phil2}). One way to view $A_b$ is to observe that
it is $\lim\limits_{\stackrel {\longleftarrow\!\!-\!\!-} 
{\empty_{p\in \mathcal{N}(A)}}}  A_p$ in $\mathsf{C^*A}$. 
While $A$ carries the subspace topology induced by the product topology
$\hspace{-6pt}\prod\limits_{p\in \mathcal{N}(A)}\hspace{-6pt} A_p$, 
$A_b$ is equipped with the
topology of ``uniform convergence," as expected. Thus,
the inclusion $(A_b,\|\cdot\|_\infty) \rightarrow A$ is continuous.
If $B \subseteq A$ is a closed $*$-algebra of $A$, then
$B_b=A_b\cap B$ (cf. \cite[3.4]{Schmu}).

\subsection{Spectrally bounded elements}
The {\em spectral radius} of an element $a$ in a pro-$C^*$-algebra 
$A$ is $r_A(a)= \sup \{ |\lambda  | \ \mid \ \lambda \in \spec_A(a)\}=
\sup\limits_{p \in \mathcal{N}(A)}  r_{A_p}(a_p)$, and $a$ is 
{\em spectrally bounded} if $r_A(a)< \nolinebreak \infty$. 
(If $A$ is not unital, $\spec_{A^+}$ is used instead of $\spec_A$.)
We put $A_{sb}$ for
the {\em set} of spectrally bounded elements in $A$ (it fails to be a
subalgebra). For every $a \in A$,  $r_A(a)\leq \|a\|_\infty$ holds, 
and so $A_b \subseteq A_{sb}$.
However, if $a$ is normal, then
$p(a)=r_{A_p}(a)$, and so $\|a\|_\infty = r_A(a)$. In particular,
a normal element is bounded if and only if it is spectrally bounded.
Therefore, $A_b$ is spanned by the spectrally bounded self-adjoint
elements in $A$ (i.e., for $A_s=\{ a\in A \mid a=a^*\}$,
$A_{sb}\cap A_s$ spans $A_b$), because each  $a \in A$ can be written 
as  $a=a_1 + i a_2$ with $a_1,a_2 \in A_s$. Since the spectrum
of an element does not depend on the topology of $A$, we obtain:

\begin{corollary} \label{cor:Ab:alg}
For a pro-$C^*$-algebra $A$, $A_b$ depends only on the algebraic structure 
of the underlying $*$-algebra of $A$. In other words, if $A_1$ and $A_2$ 
are pro-$C^*$-algebras with the same underlying $*$-algebra
(i.e., only the topologies differ), then $(A_1)_b=(A_2)_b$. \qed
\end{corollary}

\subsection{Example for $\boldsymbol{A_{sb} \neq A_b}$}
In the $\sigma$-$C^*$-algebra $A\! =\!\!\prod\limits_{n=1}^\infty 
\mathbb M_n (\mathbb{C})$ (see \ref{exm:basics}(A)), 
consider the quasi-nil\-po\-tent element $L=(L_n)$ 
given by
{\small
\begin{equation}
L_{n+1} = \left(
\begin{array}{cccccc} 
0 & 1 & 0&  \hdots  &  0& 0\\
0 & 0 & 2& 0 & .  & 0 \\
\vdots & \vdots & \ddots & \ddots & \vdots & \vdots \\
0 & 0 & \hdots& 0 & n-1 & 0\\
0 & 0 & \hdots& 0 & 0 & n\\
0 & 0 & \hdots & 0 & 0 & 0 
 \end{array} \right). 
\end{equation}}Since 
each $L_n$ is nilpotent, $\spec_{M_n (\mathbb{C})}(L_n)=\{0\}$,
and so $L$ is spectrally bounded. But, on the other hand,
\begin{equation}
n \leq \sqrt{r_{M_{n+1} (\mathbb{C})} (L_{n+1} L^*_{n+1}) }\leq 
\|L\|_\infty, 
\end{equation}
which implies that $L$ cannot be bounded. Therefore, 
$A_{sb} \neq A_b$.

\subsection{Example for $\boldsymbol {A_{sb}=A }$ } 
Let $X=[0,1]$, and put $\mathcal{D}$ for the collection of compact 
countable subsets of $X$. Since $X$ is metrizable, $\mathcal{D}$ is a
strongly functionally generating family, and thus 
$A=C_\mathcal{D}([0,1])$ is a pro-$C^*$-algebra that is not 
a $C^*$-algebra. Clearly, $A_{sb}=A_b=A$ (as sets), nevertheless, $A$ is 
not a $C^*$-algebra. This also shows that $C_\mathcal{D}([0,1])$ is not a 
$\sigma$-$C^*$-algebra:



\begin{proposition} 
If $A$ is a $\sigma$-$C^*$-algebra such that every element of $A$ is 
spectrally bounded, then $A$ is a $C^*$-algebra. 
\end{proposition}

\begin{proof}
The inclusion $(A_b,\| \cdot \|_\infty) \rightarrow A$ is always continuous. Since
$A_{sb}=A$, one has $A_b=A$, because $A_b$ is spanned by 
$A_{sb}\cap A_s$. Thus, the identity map 
$A \rightarrow (A_b,\| \cdot \|_\infty)$ 
is a $*$-homomorphism from a $\sigma$-$C^*$-algebra into a $C^*$-algebra, 
and therefore continuous by~\ref{sub:auto:sigma}.
\end{proof}

\subsection{The functor $\boldsymbol {(-)_b}$} If $\varphi\colon A 
\rightarrow B$ is a $*$-homomorphism of $*$-algebras, 
then $\varphi(A_s) \subseteq B_s$ and $\varphi(A_{sb}) \subseteq B_{sb}$, 
because $\spec_B(\varphi(a)) \subseteq \spec_A(a)$ for every $a \in A$. Thus, 
$\varphi(A_s \cap A_{sb}) \subseteq\nolinebreak B_s\cap\nolinebreak B_{sb}$. 
Therefore, if $A$ and $B$ are pro-$C^*$-algebras, then 
$\varphi(A_b) \subseteq B_b$. Hence, $\varphi$ induces a $*$-homomor\-phism
$\varphi_b\colon A_b \rightarrow B_b$. Since no assumptions were made 
about the continuity of $\varphi$, this shows that $(-)_b\colon 
\mathsf{\overline{P^*A}_d} \longrightarrow \mathsf{C^*A}$ is a functor, where  
$\mathsf{\overline{P^*A}_d}$ is the category of pro-$C^*$-algebras
and their (not necessarily continuous) $*$-homomorphisms
(cf. \cite[Folg.~3.3]{Schmu} and \cite[1.13]{Phil2}).

\begin{theorem} \label{thm:b:coref}
$(-)_b\colon \mathsf{\overline{P^*A}_d} \longrightarrow \mathsf{C^*A}$ is a 
coreflector. 
\end{theorem}

\begin{proof}
Let $B$ be a $C^*$-algebra, $A$ be a pro-$C^*$-algebra, and
$\varphi\colon B \rightarrow A$ be a  $*$-homomorphism. Then
$\varphi_b\colon B_b \rightarrow\nolinebreak A_b$ is a $*$-homomorphism,
and therefore $\varphi$ factors through $A_b$ uniquely,
because $B_b=\nolinebreak B$.
\end{proof}

\begin{corollary} \label{cor:Ab:Cb}
Let $A$ be a unital pro-$C^*$-algebra and $a \in A$ be normal.
Then $\|f(a)\|_\infty \leq \|f\|_\infty$ for every $f \in C(\spec_A(a))$.
In particular, if $f$ is bounded, then so is $f(a)$.
\qed
\end{corollary}

\subsection{Exact sequences in $\boldsymbol{\mathsf{\overline{P^*A}_d}}$}
\label{sub:exact}
A sequence
\begin{equation} \label{eq:ex:seq}
\begin{CD}
\cdots @>>> A @>\alpha>> B @>\beta>> C @>>> \cdots
\end{CD}
\end{equation}
in $\mathsf{\overline{P^*A}_d}$ is {\em exact at $B$} if
$\ker \beta = \alpha (A)$.  A categorical explanation for this choice of 
definition of exactness is that 
it is inherited from $\mathsf{T^*A_d}$, the category of
topological $*$-algebras with (not necessarily continuous) 
$*$-homomorphisms. See also Remark~\ref{rem:exact}.

\begin{theorem} \label{thm:Ab:ex}
The functor 
$(-)_b\colon \mathsf{\overline{P^*A}_d} \longrightarrow \mathsf{C^*A}$ 
preserves exactness of sequences.
\end{theorem}

\begin{proof} 
Consider the sequence in (\ref{eq:ex:seq}).
Since $\beta\alpha =0$, one has $\beta_b \alpha_b=0$, and so
$\alpha_b(A_b) \subseteq \ker \beta_b$. 
Observe that $\alpha_b(A_b)$ is closed in $B_b$, because 
$\alpha_b$ is a $*$-homomorphism between $C^*$-algebras.
Thus, in order to show the reverse inclusion, it suffices
to show that every self-adjoint $b \in \ker \beta_b$ is the 
limit of a sequence from $\alpha_b(A_b)$. Every $b \in \ker \beta_b$
has the form $b=\alpha(a)$ for some $a \in A$
because (\ref{eq:ex:seq}) is exact at $B$. If $b$ is self-adjoint, then
by replacing $a$ with $\dfrac{a+a^*} 2$ if necessary, we may assume 
that $a$ is self-adjoint too. By Corollary~\ref{cor:Ab:Cb}, $f_n(a) \in A_b$
for
\begin{equation} \label{eq:fn}
f_n(x) = \frac{n^2 x}{ n^2 + x^2}.
\end{equation}
Since $f_n (x)\rightarrow x $ uniformly on $[-\|b\|_\infty, \|b\|_\infty]$, one has 
$\alpha_b(f_n(a))=f_n(b) \stackrel {\| \cdot \|_\infty} \longrightarrow b$. 
(Observe that $f_n(a)$ is a rational polynomial in $a$, which means that 
$\alpha(f_n(a))$ is determined by $\alpha(a)$ even if $\alpha$ is not 
continuous).
Therefore, $b \in \alpha_b(A_b)$, as desired.
\end{proof}

\begin{remark}  \label{rem:exact}
The notion of exactness provided in \ref{sub:exact} might appear 
counterintuitive at first sight, because the $*$-homomorphic image 
$\alpha(A)$ need not be complete  (e.g., 
$B=C_{\mathcal{K}(\mathbb{R})}(\mathbb{R})$, $A=C_b(\mathbb{R})$,
and $\alpha=\operatorname{id}$). Requiring 
$\overline{\alpha(A)}=\ker \beta$ instead could make the impression of 
being a better definition, however, this is not the case:
For $A=C([0,1])$, $B=\mathbb{C}^{[0,1]}$, and $\alpha$ the inclusion,
$\overline{\alpha(A)}$ is dense in $B$, but $\alpha_b(A)$ is a closed
proper subalgebra of $B_b=l^\infty([0,1])$. Therefore, 
$\alpha(A)=\ker \beta$ is indeed necessary in 
Theorem~\ref{thm:Ab:ex}.
\end{remark}

\begin{corollary}
Let $I$ be a closed $*$-ideal of the pro-$C^*$-algebra $A$ such that
$A/I$ is complete (i.e., $A/I$ is a pro-$C^*$-algebra).
Then $I_b$ is a closed $*$-ideal of $A_b$, and $(A/I)_b \cong A_b/I_b$.
\end{corollary}

\begin{proof}
Apply Theorem~\ref{thm:Ab:ex} to the sequence
$0 \rightarrow I \rightarrow A \rightarrow A/I \rightarrow 0$.
\end{proof}

\subsection{The quotient $\boldsymbol{A/I}$ need not be complete}
Because of the automatic continuity property of $\sigma$-$C^*$-algebras
(see \ref{sub:auto:sigma}), every short exact sequence of
$\sigma$-$C^*$-algebras is the limit of a countable collection
of short exact sequences of $C^*$-algebras (cf. \cite[5.3]{Phil1}).
In particular, if $A$ is a $\sigma$-$C^*$-algebra and $I$ is its closed 
$*$-ideal, then $A/I$ is complete (cf. \cite[5.4]{Phil1}).
The quotient of a complete topological algebra by a closed ideal
need not be complete (cf. \cite{Koethe}), but one would hope that 
pro-$C^*$-algebras are better behaved. This is, however, not the case, as 
the following example shows: Let $X$ be a locally compact space that is 
not normal. Then $X$ contains a closed subspace $Y \subseteq X$ such that 
not every continuous map on $Y$ can be extended continuously to $X$. Put
$A=C_{\mathcal{K}(X)}(X)$ and $B=C_{\mathcal{K}(Y)}(Y)$, let
$\varphi\colon A \rightarrow B$ be  the restriction (i.e., 
$f \mapsto f \rest_Y$), and set $I=\ker \varphi$. Since $Y$ is closed in 
$X$, it is easy to see that $\varphi$ is open onto its image, and thus
$A/I \cong \varphi(A)$ as topological $*$-algebras. 
The image, $\varphi(A)$, is precisely the $*$-subalgebra of continuous 
maps on $Y$ that can be continuously extended to $X$, and so 
$\varphi(A)\neq B$. On the other hand, by the  Stone-Weierstrass Theorem, 
$\varphi(A)$ is dense in $B$. Therefore, $\varphi(A)\cong A/I$ is 
not complete. (This example, which is based on   
\cite[Exercise~IV.11]{SchaefWol}, was  suggested by Subhash Bhatt, and was 
communicated personally to the authors by N. C. Phillips.)


\begin{corollary} \label{cor:Ab:Ap}
Let $A$ be a pro-$C^*$-algebra. Then
$A_p\cong A_b/(\ker p)_b$ for every $p \in \mathcal{N}(A)$.
\qed
\end{corollary}

\begin{proposition}
Let $A$ be a pro-$C^*$-algebra. If $A_b$ is simple, then $A$ is a 
$C^*$-algebra.
\end{proposition}

\begin{proof}
Since $A_b$ is dense in $A$ (see \ref{sub:Ab}), 
$(\ker p)_b \neq A_b$ for every nonzero $p \in \mathcal{N}(A)$.
Thus, $(\ker p)_b = 0$, because $A_b$ is simple, and 
and therefore $\ker p =0$, because $(\ker p)_b$ is dense in $\ker p$.
Hence, by Corollary~\ref{cor:Ab:Ap}, $A=A_p = A_b$
for every nonzero $p \in \mathcal{N}(A)$, which 
means that $A$ is a $C^*$-algebra.
\end{proof}

\subsection{Problem} Let $B$ be a unital $C^*$-algebra that is not simple.
Is there a pro-$C^*$-algebra $A$ that is not a $C^*$-algebra
such that $B=A_b$? This problem is fairly non-trivial even in the abelian 
case, where $B=C(K)$ for some compact Hausdorff space $K$. If $K$ is 
metrizable and uncountable, then $A=C_{\mathcal{D}}(K)$ is a positive 
solution, where $\mathcal{D}$ is the collection of compact countable 
subsets of $K$. In fact, it is sufficient to assume that $K$ is 
compact and a {\em sequential space} (i.e., continuity of maps is 
determined by continuity on sequences). Nevertheless, we do not have 
an answer for the problem even for such well-known $C^*$-algebras as
$c$ and $c_0$ (the subalgebras of $l^\infty$ consisting of the convergent 
and zero sequences, respectively).

\section{The group of unitary elements}
\label{sect:u}

In this section, we are concerned with the topological group $A_u$ of
unitary elements in a unital pro-$C^*$-algebra $A$ (i.e., 
$x x^* =x^* x = 1$), and study the connected component $A_{u,1}$ of $1$
in it. A key tool in our investigation is the subgroup 
$A_{\exp}  = \{ e^{ia_1}\cdots e^{i a_n} \mid a_1,\ldots,a_n \in A_s\}$, 
generated by exponentials of self-adjoint elements. 

\begin{proposition} \label{prop:unit:1}
Let $A$ be a pro-$C^*$-algebra. Then:

\begin{list}{{\rm (\alph{enumi})}}
{\usecounter{enumi}\setlength{\labelwidth}{25pt}\setlength{\topsep}{-8pt}
\setlength{\itemsep}{-5pt} \setlength{\leftmargin}{25pt}}

\item
$A_{\exp}$ is path-connected in the topology of $A$;

\item
$A_{\exp} \subseteq A_{u,1}$;

\item
$A_{\exp}$ is an clopen subgroup of $A_u$ in the 
$\|\cdot\|_\infty$-topology.

\end{list}
\end{proposition}

\begin{proof}
We note that since a unitary element has norm $1$ in any $C^*$-algebra,
$\|u\|_\infty=1$ for every $u \in A_u$. So $A_u \subseteq A_b$, and 
$A_u=(A_b)_u$ as sets (their topology might be quite different, though).

(a) We show that each generator $e^{ia}$ ($a \in A_s$) can
be joined to $1$ by a path. To that end, consider the 
map $u\colon [0,1]\rightarrow A$ given by
$u(t)=e^{ita}$. It suffices to show that $u$ is continuous. 
For each $p \in \mathcal{N}(A)$,
$f(t,s)=e^{its}$ is uniformly continuous on the compact set
$[0,1] \times \spec_{A_p}(a_p)$, and thus 
$u_p(t) = e^{ita_p}$ is continuous $[0,1]\rightarrow A_p$.
Therefore, $u$ is continuous, as desired.

(b) follows from (a).

(c) Since every open subgroup in a topological group is also
closed, it suffices to show that 
$\{ u \in A_u \mid \|1 -u \|_\infty < 1\} \subseteq A_{\exp}$.
Let $u \in A_u$ be such that $\|1-u \|_\infty < 1$. Then
$\spec_A(u) \subseteq \mathbb{S}\backslash \{-1\}$ (where 
$\mathbb{S}\! =\! \{z \in \mathbb{C}\! \mid |z|=1\}$), and so
$f(z)=\arg z$ can be continuously defined on $\spec_A(u)$.
Thus, by the functional calculus, $u=e^{ia}$
for $a=f(u)$, and $a \in A_s$ because $f$ is real-valued.
Therefore, $u \in A_{\exp}$, as desired.
\end{proof}

\begin{remark} \label{rem:C*:exp=1}
If $A$ is a $C^*$-algebra, then $A_{\exp}$ is a  clopen 
path-connected subgroup of $A_u$ (by Proposition~\ref{prop:unit:1}(c)), 
and hence $A_{\exp}=A_{u,1}$.
\end{remark}

\begin{lemma}
Let $\alpha\colon A \rightarrow B$ be continuous $*$-homomorphism
of pro-$C^*$-algebras. Then  one has $\alpha(g(a))=g(\alpha(a))$
for every normal $a \in A$ and $g \in C(\spec_A(a))$.
\end{lemma}

\begin{proof}
Put $b=\alpha(a)$. For every polynomial $p(z,\bar z)$,
$\alpha(p(a,a^*))=p(b,b^*)$, because $\alpha$ is a $*$-homomorphism.
If $g \in C(\spec_A(a))$, then it can be uniformly approximated
by polynomials $p_\alpha(z,\bar z)$ on compact sets, and therefore
$g(a)=\lim p_\alpha(a,a^*)$ and $g(b)=\lim p_\alpha(b,b^*)$.
Hence, the statement follows by continuity of $\alpha$.
\end{proof}

\begin{corollary} \label{cor:unit:onto}
Let $\alpha\colon A \rightarrow B$ be a surjective continuous $*$-homomorphism
of pro-$C^*$-algebras. Then 
$\alpha (A_{\exp})=B_{\exp}$. \qed
\end{corollary}

\begin{theorem}
Let $A$ be a pro-$C^*$-algebra. Then:

\begin{list}{{\rm (\alph{enumi})}}
{\usecounter{enumi}\setlength{\labelwidth}{25pt}\setlength{\topsep}{-8pt}
\setlength{\itemsep}{-5pt} \setlength{\leftmargin}{25pt}}

\item
$A_{\exp}$ is dense in $A_{u,1}$;

\item
$A_{u,1}=\lim\limits_{\stackrel {\longleftarrow\!\!-\!\!-} 
{\empty_{p\in \mathcal{N}(A)}}} (A_p)_{u,1}$.

\end{list}
\end{theorem}

\begin{proof}
To shorten notations, put $G=\!\!
\lim\limits_{\stackrel {\longleftarrow\!\!-\!\!-} 
{\empty_{p\in \mathcal{N}(A)}}} (A_p)_{u,1}$.
Clearly,  $A_{u,1} \subseteq G \subseteq A_u$.

(a) We show that $A_{\exp}$ is dense in $G$. Since each projection 
$\pi_p\colon A \rightarrow A_p$ is surjective, by 
Corollary~\ref{cor:unit:onto}, $\pi_p(A_{\exp})=(A_p)_{\exp}$.
Since $A_p$ is a $C^*$-algebra, $(A_p)_{\exp}=(A_p)_{u,1}$
(cf. Remark~\ref{rem:C*:exp=1}), and thus $\pi_p(A_{\exp})=(A_p)_{u,1}$.
Therefore, the statement follows, because $\mathcal {N}(A)$ is directed.

(b) Since $A_{\exp}\subseteq A_{u,1}$, $A_{u,1}$ is dense in $G$.
On the other hand, $A_{u,1}$ is a closed subgroup of $A_u$,
because it is a connected component. Therefore, $A_{u,1}=G$,
as desired.
\end{proof}

\subsection{Example for $\boldsymbol{A_{\exp}\neq A_{u,1}}$}
Let $B=C([0,1]^\omega)\otimes M_2(\mathbb{C})$ be the $C^*$-algebra of
continuous maps $f\colon [0,1]^\omega \rightarrow M_2(\mathbb{C})$ from 
the Hilbert cube into $M_2(\mathbb{C})$ (or equivalently, of $2 \times 2$ 
matrices with entries in $C([0,1]^\omega)$). Phillips showed that
for every $n \in \mathbb{N}$, there is $u_n\in B_u$ in the {\em path} 
component of $1$ in $B_u$ such that $u_n$ cannot be expressed as the 
product of $n$ many exponentials $e^{ib}$, where $b \in B_s$
(cf.~\cite[2.6]{Phil4}). In other 
words, $u_n\neq e^{ib_1}\cdots e^{ib_n}$ for every $b_1,\ldots,b_n \in B_s$.
Put $A=B^\omega$, the countable power of $B$ in $\mathsf{T^*A}$;
$A$ is a $\sigma$-$C^*$-algebra, and $u=(u_n) \in A_u$.
If $p_n\colon \mathbb{I} \rightarrow B_u$ is a path such that
$p_n(0)=1$ and $p_n(1)=u_n$, then $p=(p_n)\colon \mathbb{I} \rightarrow A_u$ 
is a path joining $u$ to $1$ in $A_u$. Thus, $u$ is in the path component 
of $A_u$, but it cannot be in $A_{\exp}$, because $u$ is not a product of 
$n$ exponentials for any $n\in \mathbb{N}$. In particular, 
$A_{\exp} \subsetneq A_{u,1}$.
(This example is due to N. C. Phillips \cite[2.11]{Phil4}.)

\section*{Acknowledgments}

We wish to thank Professor  N. C. Phillips for the valuable 
discussions that were of great assistance in writing this paper.

{\footnotesize

\bibliography{proCalg}

\def\cprime{$'$}
\begin{thebibliography}{10}

\bibitem{Allan}
Graham~R. Allan.
\newblock Stable inverse-limit sequences, with application to {F}r\'echet
  algebras.
\newblock {\em Studia Math.}, 121(3):277--308, 1996.

\bibitem{Apost}
Constantin Apostol.
\newblock {$b\sp{\ast} $}-algebras and their representation.
\newblock {\em J. London Math. Soc. (2)}, 3:30--38, 1971.

\bibitem{ArhBook}
A.~V. Arkhangel{\cprime}ski{\u\i}.
\newblock {\em Topological function spaces}, volume~78 of {\em Mathematics and
  its Applications (Soviet Series)}.
\newblock Kluwer Academic Publishers Group, Dordrecht, 1992.
\newblock Translated from the Russian by R. A. M. Hoksbergen.

\bibitem{Arveson}
William Arveson.
\newblock The harmonic analysis of automorphism groups.
\newblock In {\em Operator algebras and applications, Part I (Kingston, Ont.,
  1980)}, volume~38 of {\em Proc. Sympos. Pure Math.}, pages 199--269. Amer.
  Math. Soc., Providence, R.I., 1982.

\bibitem{BhaKar}
Subhash~J. Bhatt and Dinesh~J. Karia.
\newblock An intrinsic characterization of {P}ro-{$C\sp *$}-algebras and its
  applications.
\newblock {\em J. Math. Anal. Appl.}, 175(1):68--80, 1993.

\bibitem{Board}
Michael~E. Boardman.
\newblock Relative spectra in complete lmc-algebras with applications.
\newblock {\em Illinois J. Math.}, 39(1):119--139, 1995.

\bibitem{Bourd2}
G.~Bourdaud.
\newblock Quelques aspects de la dualit\'e en analyse fonctionnelle.
\newblock {\em Diagrammes}, 5:B1--B12, 1981.

\bibitem{Bourd1}
G{\'e}rard Bourdaud.
\newblock Sur la dualit\'e des alg\`ebres localement convexes.
\newblock {\em C. R. Acad. Sci. Paris S\'er. A-B}, 281(23):Ai, A1011--A1014,
  1975.

\bibitem{Brown}
R.~Brown.
\newblock Function spaces and product topologies.
\newblock {\em Quart. J. Math. Oxford Ser. (2)}, 15:238--250, 1964.

\bibitem{EHR1Chi}
M.~Chidami and R.~El~Harti.
\newblock Calcul fonctionnel holomorphe en dimension infinie dans les lmca.
\newblock {\em Rend. Circ. Mat. Palermo (2)}, 48(3):541--548, 1999.

\bibitem{DixmierC}
Jacques Dixmier.
\newblock {\em {$C\sp*$}-algebras}.
\newblock North-Holland Publishing Co., Amsterdam, 1977.
\newblock Translated from the French by Francis Jellett, North-Holland
  Mathematical Library, Vol. 15.

\bibitem{DubPor}
Eduardo~J. Dubuc and Horacio Porta.
\newblock Convenient categories of topological algebras, and their duality
  theory.
\newblock {\em J. Pure Appl. Algebra}, 1(3):281--316, 1971.

\bibitem{DubPor2}
Eduardo~J. Dubuc and Horacio Porta.
\newblock Uniform spaces, {S}panier quasitopologies, and a duality for locally
  convex algebras.
\newblock {\em J. Austral. Math. Soc. Ser. A}, 29(1):99--128, 1980.

\bibitem{EHR4}
Rachid El~Harti.
\newblock Contractible {F}r\'echet algebras.
\newblock {\em Proc. Amer. Math. Soc.}, 132(5):1251--1255 (electronic), 2004.

\bibitem{Engel6}
Ryszard Engelking.
\newblock {\em General topology}, volume~6 of {\em Sigma Series in Pure
  Mathematics}.
\newblock Heldermann Verlag, Berlin, second edition, 1989.
\newblock Translated from the Polish by the author.

\bibitem{Inoue}
Atsushi Inoue.
\newblock Locally {$C\sp{\ast} $}-algebra.
\newblock {\em Mem. Fac. Sci. Kyushu Univ. Ser. A}, 25:197--235, 1971.

\bibitem{Kelley}
John~L. Kelley.
\newblock {\em General topology}.
\newblock D. Van Nostrand Company, Inc., Toronto-New York-London, 1955.

\bibitem{Koethe}
Gottfried K{\"o}the.
\newblock {\em Topological vector spaces. {I}}.
\newblock Translated from the German by D. J. H. Garling. Die Grundlehren der
  mathematischen Wissenschaften, Band 159. Springer-Verlag New York Inc., New
  York, 1969.

\bibitem{GL5}
G{\'a}bor Luk{\'a}cs.
\newblock A convenient subcategory of {T}ych.
\newblock {\em Appl. Categ. Structures}, 12(4):369--377, 2004.

\bibitem{GL7}
G\'abor Luk\'acs.
\newblock Lifted closure operators.
\newblock {\em Preprint}, ArXiv: math.CT/0502410.

\bibitem{MacLane}
Saunders Mac~Lane.
\newblock {\em Categories for the working mathematician}, volume~5 of {\em
  Graduate Texts in Mathematics}.
\newblock Springer-Verlag, New York, second edition, 1998.

\bibitem{Michael}
Ernest~A. Michael.
\newblock Locally multiplicatively-convex topological algebras.
\newblock {\em Mem. Amer. Math. Soc.,}, 1952(11):79, 1952.

\bibitem{Palmer1}
Theodore~W. Palmer.
\newblock {\em Banach algebras and the general theory of {$\sp *$}-algebras.
  {V}ol. {I}}, volume~49 of {\em Encyclopedia of Mathematics and its
  Applications}.
\newblock Cambridge University Press, Cambridge, 1994.
\newblock Algebras and Banach algebras.

\bibitem{Palmer2}
Theodore~W. Palmer.
\newblock {\em Banach algebras and the general theory of {$*$}-algebras. {V}ol.
  2}, volume~79 of {\em Encyclopedia of Mathematics and its Applications}.
\newblock Cambridge University Press, Cambridge, 2001.
\newblock $*$-algebras.

\bibitem{PelRos2}
J.~Wick Pelletier and J.~Rosick{\'y}.
\newblock On the equational theory of {$C\sp *$}-algebras.
\newblock {\em Algebra Universalis}, 30(2):275--284, 1993.

\bibitem{PelRos}
Joan~Wick Pelletier and Ji{\v{r}}{\'\i} Rosick{\'y}.
\newblock Generating the equational theory of {$C\sp *$}-algebras and related
  categories.
\newblock In {\em Categorical topology and its relation to analysis, algebra
  and combinatorics (Prague, 1988)}, pages 163--180. World Sci. Publishing,
  Teaneck, NJ, 1989.

\bibitem{Phil2}
N.~Christopher Phillips.
\newblock Inverse limits of {$C\sp *$}-algebras.
\newblock {\em J. Operator Theory}, 19(1):159--195, 1988.

\bibitem{Phil1}
N.~Christopher Phillips.
\newblock Inverse limits of {$C\sp *$}-algebras and applications.
\newblock In {\em Operator algebras and applications, Vol.\ 1}, volume 135 of
  {\em London Math. Soc. Lecture Note Ser.}, pages 127--185. Cambridge Univ.
  Press, Cambridge, 1988.

\bibitem{Phil3}
N.~Christopher Phillips.
\newblock Representable {$K$}-theory for {$\sigma$}-{$C\sp *$}-algebras.
\newblock {\em $K$-Theory}, 3(5):441--478, 1989.

\bibitem{Phil4}
N.~Christopher Phillips.
\newblock How many exponentials?
\newblock {\em Amer. J. Math.}, 116(6):1513--1543, 1994.

\bibitem{SchaefWol}
H.~H. Schaefer and M.~P. Wolff.
\newblock {\em Topological vector spaces}, volume~3 of {\em Graduate Texts in
  Mathematics}.
\newblock Springer-Verlag, New York, second edition, 1999.

\bibitem{Schmu}
Konrad Schm{\"u}dgen.
\newblock \"{U}ber {LMC}-{A}lgebren.
\newblock {\em Math. Nachr.}, 68:167--182, 1975.

\bibitem{Sebes}
Z.~Sebesty{\'e}n.
\newblock Every {$C\sp{\ast} $}-seminorm is automatically submultiplicative.
\newblock {\em Period. Math. Hungar.}, 10(1):1--8, 1979.

\bibitem{Steenrod}
N.~E. Steenrod.
\newblock A convenient category of topological spaces.
\newblock {\em Michigan Math. J.}, 14:133--152, 1967.

\bibitem{Osdol}
D.~H. Van~Osdol.
\newblock {$C\sp \ast$}-algebras and cohomology.
\newblock In {\em Categorical topology (Toledo, Ohio, 1983)}, volume~5 of {\em
  Sigma Ser. Pure Math.}, pages 582--587. Heldermann, Berlin, 1984.

\bibitem{Zelaz}
Wies{\l}aw {\.Z}elazko.
\newblock {\em Banach algebras}.
\newblock Elsevier Publishing Co., Amsterdam-London-New York, 1973.
\newblock Translated from the Polish by Marcin E. Kuczma.

\end{thebibliography}
}

\bigskip
\noindent
\begin{tabular}{l @{\hspace{1.9cm}} l}
Department of Mathematics & Department of Mathematics and Statistics\\
University Hassan I, FST de Settat & Dalhousie University \\
BP 577, 2600 Settat & Halifax, B3H 3J5, Nova Scotia \\
Morocco & Canada \\ & \\
\em e-mail: elharti@ibnsina.uh1.ac.ma  &
\em e-mail: lukacs@mathstat.dal.ca
\end{tabular}

\end{document}